\documentclass[oneside]{amsart}
\usepackage{a4wide}
\usepackage{graphicx}
\usepackage{amssymb}

\makeatletter
 \theoremstyle{plain}    
 \newtheorem{thm}{Theorem}[section]
 \numberwithin{equation}{section} 
 \numberwithin{figure}{section} 
 \theoremstyle{plain}
 \theoremstyle{plain}    
 \newtheorem{lem}[thm]{Lemma} 
 \theoremstyle{plain}    
 \newtheorem{prop}[thm]{Proposition} 

\AtBeginDocument{
  
}

\makeatother
\begin{document}

\title{On deformation quantization of Dirac structures}
\author{Pavol \v Severa}
\address{Section de Math\'ematiques, rue du Li\`evre 2-4, 1211 Geneva,
Switzerland\\ on leave from Dept.~of Theor.~Physics, Mlynsk\'a dolina F2,
84248 Bratislava, Slovakia}
\email{pavol.severa@gmail.com}
\thanks{supported in part by the Swiss National Science Foundation}

\begin{abstract}
Motivated by the problem of transverse deformation quantization of  foliated
manifolds, we describe a quantization of Dirac structures (more precisely,
of those that are formal deformations
of regular ones) to stacks of algebroids in the sense of Kontsevich.

\end{abstract}

\maketitle

\section{Introduction}

This paper is, roughly speaking, about transverse deformation quantization.
If $N$ is a foliated manifold, we have the sheaf of the functions
constant along the leaves. If we suppose that on $N$ we also have
a transverse Poisson structure (that is, we make the above-mentioned
sheaf to a sheaf of Poisson algebras), we can locally quantize the
algebras of functions constant along the leaves; the result will,
in general, be just a {}``sheaf up to homotopy'', namely a stack
of algebras. It is then an overkill to start with a transverse Poisson
structure; it would be enough to have an {}``up to homotopy transverse
Poisson structure'', and we would still get a stack after quantization.
The correct notion of {}``up to homotopy transverse Poisson structure''
is in this context the following: a formal family of Dirac structures
$C(\hbar)$ on $N$ such that the leaves of $C(0)$ are precisely
the leaves of the foliation. This paper is thus about quantization
of Dirac structures, more precisely, about quantization of Dirac structures
that are formal deformations of regular ones.

Let thus $C(\hbar)$ be a family of Dirac structures on $N$, formally
depending on $\hbar$, with the property that $C(0)$ is regular,
i.e.~that it gives a foliation $F$ of $N$. If $M\subset N$ is
a local transversal (that is, a submanifold of dimension complementary
to $F$, intersecting the leaves transversally), the Dirac structure
$C(\hbar)$ induces on $M$ a Poisson structure of the form \[
\pi=\pi^{(1)}\hbar+\pi^{(2)}\hbar^{2}+\pi^{(3)}\hbar^{3}+\dots\]
 These local Poisson structures are glued by $C(\hbar)$ in such a
way that their quantizations form a stack on $N$. It is proved by
a straightforward use of the formality theorem of Kontsevich; this
is what we do in this paper. 

The stack (sheaf of categories) that we get is a deformation of the
following: for any open $U\subset N$ the objects are line bundles
over $U$, endowed with connection along the leaves of $F$, with
curvature equal to $C(0)$; morphisms are maps of bundles preserving
the connections. Locally (i.e.~for nice $U$'s) this category has
one object up to isomorphisms and the algebra of its endomorphisms
is the algebra of functions constant along the leaves. 

This paper leaves at least two open questions. First is: can one quantize
Dirac structures that are not deformations of regular ones? And the
second (more interesting), connected with non-commutative geometry:
Instead of using the sheaf of functions constant along the leaves,
one can form a global object, a non-commutative algebra that is supposed
to replace functions on the possibly degenerate space of leaves. What
happens with this algebra after quantization of $C(\hbar)$?

This paper is entirely based on ideas from \cite{K2}. Dirac structures
were used for a deformation quantization problem, using the same techniques
as here, in \cite{moi}. $C^*$-quantization of tori with constant Dirac
structures was introduced by Tang and Weinstein in \cite{TW}; in their simple
case they didn't have to use homotopy techniques of the type that appear
here, and they were able to get a global non-commutative algebra.

\section{Dirac structures}

Dirac structures were introduced by Courant \cite{Cou} as a common
generalization of closed 2-forms and Poisson structures. First idea
is to see 2-forms on a manifold $N$ as maps $TN\rightarrow T^{*}N$, bivectors
as maps $T^{*}N\rightarrow TN$, and represent both by their graphs, which
are subbundles in $(T\oplus T^{*})N$. Skew-symmetry of 2-forms/bivectors
makes these graphs isotropic with respect to the inner product \[
\langle(u,\alpha),(v,\beta)\rangle=\alpha(v)+\beta(u).\]
A common generalization of both bivectors and 2-forms are then maximally
isotropic subbundles of $(T\oplus T^{*})N$.

It then turns out that the conditions on a 2-form to be closed, and
on a bivector to be Poisson, can be rephrased as a condition on the
subbundle. We shall do it the standard mysterious way: one introduces
so-called Courant bracket on sections of $(T\oplus T^{*})N$, \[
[\![(u,\alpha),(v,\beta)]\!]=[u,v]+\mathcal{L}_{u}\beta-i_{v}d\alpha,\]
and then defines \emph{Dirac structures} to be maximally isotropic
subbundles of $(T\oplus T^{*})N$, closed under Courant bracket.%
\footnote{If it makes it less mysterious -- sections of $(T\oplus T^*)N$
are derivations of degree $-1$ of $\Omega(N)[t]$, where $t$ is an 
auxiliary variable of degree 2, $dt=0$; for two such derivations $\xi_{1}$, $\xi_2$
we have $[\![\xi_1,\xi_2]\!]=[[d,\xi_1],\xi_2]$ and
$\langle \xi_1,\xi_2\rangle\,\partial_t=[\xi_1,\xi_2]$.}

Similar to symplectic leaves of Poisson structures, any Dirac structure
gives us a singular foliation on $N$, with closed 2-forms on its
leaves (contrary to Poisson structures, these closed 2-forms don't
have to be symplectic). A particular case of a Dirac structure is
thus a (non-singular) foliation, with a closed 2-form on its leaves.
Such Dirac structures (i.e.~those, whose leaves have constant dimension)
are called \emph{regular}.

\section{Dirac structures as Hamiltonian families of Poisson structures}

Let $L=\Gamma(\bigwedge TM)[1]$ be the DGLA of multivector fields
on a manifold $M$ (the bracket in $L$ is the Schouten bracket and
the differential is zero). Let us take another manifold $B$ and try
to find the solutions of the MC equation in the DGLA $L\otimes\Omega(B)$.%
\footnote{Caution: here, and in similar situations, we don't mean the algebraic
tensor product, but rather its natural completion: the coefficients
of the objects from $L\otimes\Omega(B)$ are allowed to be any smooth
functions on $M\times B$.%
} We thus look for the sections $\sigma$ of degree $2$ of the vector
bundle $\bigwedge TM\otimes\bigwedge T^{*}B=\bigwedge(TM\oplus T^{*}B)$
over $M\times B$, satisfying the MC equation \[
d\sigma+[\sigma,\sigma]/2=0.\]
Such a section can be viewed as a skew-symmetric bilinear form on
$T^{*}M\oplus TB$, and thus as a maximal isotropic subbundle of $(T\oplus T^{*})(M\times B)$.
A simple computation shows that this subbundle is a Dirac structure
iff $\sigma$ solves the MC equation. We thus have the following:

\begin{lem}
A solution of the MC equation in $\Gamma(\bigwedge TM)[1]\otimes\Omega(B)$
is equivalent to a Dirac structure on $M\times B$ transversal to
$TM\oplus T^{*}B$.
\end{lem}
Now suppose we replace $L=\Gamma(\bigwedge TM)[1]$ by multivectors
formally depending on $\hbar$, as they appear in deformation quantization.
That is, let the DGLA $L'$ be given by\[
L'^{-1}=C^{\infty}(M)[[\hbar]],\]

\[
L'^{i}=\hbar\,\Gamma({\textstyle \bigwedge}^{i+1}TM)[[\hbar]]\textrm{ for }i\geq0.\]
We have the following version of the previous lemma:

\begin{lem}
\label{lem-form-dir}A solution of the MC equation in $L'\otimes\Omega(B)$
is equivalent to a (formal) family of Dirac structures $C(\hbar)$
on $M\times B$ such that $C(0)$ is a Dirac structure with leaves
$\{ x\}\times B$, $x\in M$.
\end{lem}
In analogy with Hamiltonian families of symplectic structures, we
will call solutions of the MC equation in $L'\otimes\Omega(B)$ \emph{Hamiltonian
families} of formal Poisson structures on $M$ parametrized by $B$
(I called them {}``tight families'' in \cite{moi}, but it's better
to follow an older tradition)\emph{.}

\section{Quantization of Hamiltonian families of formal Poisson structures\label{sec:Quantization-of-formal}}

Let us denote $L''$ the DGLA of polydifferential operators on $M$
formally depending on $\hbar$, more preciselly\[
L''^{-1}=C^{\infty}(M)[[\hbar]],\]

\[
L''^{i}=\hbar\,\mathit{PD}^{i+1}(M)[[\hbar]]\textrm{ for }i\geq0,\]
where $\mathit{PD}^{i}$ denotes the space of polydifferential operators
\[
\underbrace{C^{\infty}(M)\times\dots\times C^{\infty}(M)}_{i\textrm{ times}}\rightarrow C^{\infty}(M).\]

By the formality theorem, for any solution $\sigma$ of the MC equation
in $L'\otimes\Omega(B)$ we have a solution $\tau$ in $L''\otimes\Omega(B)$;
we call the latter a \emph{tight family of $*$-products} on $M$
parametrized by $B$\emph{.} 

Let us try to understand what tight $*$-product families actually
are. The rest of this section is completely stolen from \cite{K2}.
Let us decompose $\tau$ to bihomogeneous components $\tau=\tau^{0}\oplus\tau^{1}\oplus\tau^{2}$
(the superscript means the degree in $\Omega(B)$) to see what the
MC equation means; $m:C^{\infty}(M)\times C^{\infty}(M)\rightarrow C^{\infty}(M)$
denotes the ordinary product of functions:

\[
[\tau^{0}+m,\tau^{0}+m]/2=0,\]
i.e. $m+\tau^{0}$ is a family of $*$-products on $M$ parametrized
by $B$;

\[
d\tau^{0}+[\tau^{0}+m,\tau^{1}]=0,\]
i.e.~$\tau^{1}\in\Omega^{1}(B)\otimes\hbar\, D(M)[[\hbar]]$, understood
as a connection on the infinite-dimensional vector bundle $C^{\infty}(M)[[\hbar]]\times B\rightarrow B$,
makes the family of $*$-products parallel;

\[
d\tau^{1}+[\tau^{1},\tau^{1}]/2+[\tau^{0}+m,\tau^{2}]=0,\]
i.e. the curvature of $\tau^{1}$ is an inner derivation of the algebra
$C^{\infty}(M)[[\hbar]]$ with its $*$-product $\tau^{0}+m$, and
finally\[
d\tau^{2}+[\tau^{1},\tau^{2}]=0.\]

For any point $b\in B$, let us denote the algebra $C^{\infty}(M)[[\hbar]]$
with the $*$-product $m+\tau^{0}(b)$ by $A_{b}$.

Although $\tau^{1}$ is a connection on an infinite-dimensional vector
bundle, its parallel transport is well-defined (this is because $\tau^{1}\in\Omega^{1}(B)\otimes\hbar\, D(M)[[\hbar]]$,
i.e.~$\tau^{1}=O(\hbar)$). One can also use $\tau^{2}$ to get a
{}``2-dimensional parallel transport'', and we get the following:

\begin{enumerate}
\item for every curve $\gamma$ in $B$, connecting points $b_{1}$ and
$b_{2}$, an isomorphism $T_{\gamma}:A_{b_{1}}\rightarrow A_{b_{2}}$; $T_{\gamma}$
is just the parallel transport along $\gamma$,
\item for every disk $D$ in $B$ with a chosen point $b$ on the boundary,
an element $a_{D,b}\in A_{b}$.
\end{enumerate}
The following relations then hold:

\begin{enumerate}
\item if $\gamma$ is the boundary of $D$ (so that $b_{1}=b_{2}=b$) then
$T_{\gamma}$ is the inner automorphism given by $a_{D,b}\in A_{b}$
\item $a_{D,b}\in A_{b}$ depends only on the homotopy class of $D$ rel
boundary
\item $a_{D,b}$'s are multiplicative ($a_{D_{1}\cup D_{2},b}=a_{D_{1},b}\, a_{D_{2},b}$)
and behave naturally under change of $b$ ($a_{D,b_{2}}=T_{\gamma}\  a_{D,b_{1}}$
where $\gamma$ is the curve from $b_{1}$ to $b_{2}$), see the picture:
\end{enumerate}
\begin{center}\includegraphics[%
  scale=0.5]{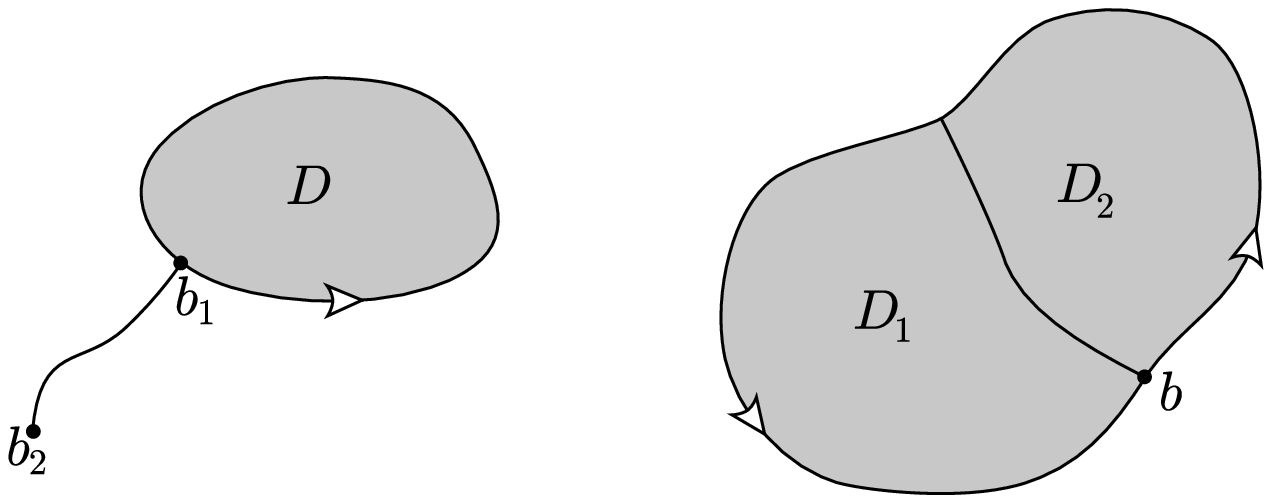}\end{center}

The quantization of Hamiltonian families to tight $*$-product families,
as we have just described it, depends on the choice of the quasiisomorphism
$L'\rightarrow L''$. Unfortunately, there is no natural (i.e.~diffeomorphisms-invariant)
choice for the quasiisomorphism, but fortunately, it is natural up
to homotopy. Using the precise meaning of the previous sentence, given
by Kontsevich in \cite{K2}, we have:

\begin{prop}
\label{pro:connection}Given a Hamiltonian family of formal Poisson
structures on $M$ parametrized by $B$, and given also a smooth family
of connections on $M$ parametrized by $B$, there is a natural and
local construction of a tight family of $*$-products on $M$ parametrized
by $B$.
\end{prop}

\section{Deformation quantization of Dirac structures}

Let now $N$ be a manifold and $C(0)$ a regular Dirac structure on
$N$. It means that we are given a foliation $F$ of $N$ and a closed
2-form on the leaves of $F$. Suppose we extend $C(0)$ to a formal
family $C(\hbar)$ of Dirac structures.

Let us choose a connection on the vector bundle $TN/F$ (the normal
bundle of the foliation), to be able to use Proposition \ref{pro:connection}.
If $M\subset N$ is a local transversal, i.e.~a submanifold with
complementary dimension to $F$ and transversal to the leaves, the
connection induces a connection on $TM$. The Dirac structure $C(\hbar)$
pulls back to a Poisson structure on $M$ of the form \[
\pi=\pi^{(1)}\hbar+\pi^{(2)}\hbar^{2}+\pi^{(3)}\hbar^{3}+\dots\]
 More generally, for any smooth map $f:M\times B\rightarrow N$ such that 

\begin{itemize}
\item for any point $p=(m,b)\in M\times B$, $d_{p}f$ maps $T_{b}B$ into
$F_{f(p)}$, and the induced map $T_{m}M\rightarrow(T_{f(p)}N)/F_{f(p)}$
is a bijection,
\end{itemize}
the pullback of $C(\hbar)$ is a formal family of Dirac structures
on $M\times B$ satisfying the conditions of Lemma \ref{lem-form-dir}.
Together with the pullback of the connection, it gives us (via Proposition
\ref{pro:connection}) a tight $*$-product family on $M$ parametrized
by $B$, and thus the isomorphisms $T$ and elements $a$ as is Section \ref{sec:Quantization-of-formal}.

From here it is not difficult to see that the quantizations of the
Poisson structures on transversals form a stack on $N$. One needs
the cases $B=I$ ($I$ is an interval), i.e.~homotopies of transversals,
$B=I\times I$, i.e.~homotopies of homotopies, and finally $B=I\times I\times I$.
We shall construct the stack in some detail in what follows.

\subsection{Algebroids}

A \emph{linear category} is a category in which $\mathrm{Hom}(X,Y)$
is a vector space (for any two objects $X$, $Y$) and composition
of morphisms is bilinear. An \emph{algebroid} is a linear category
in which any two objects are isomorphic. Any algebroid gives us an
algebra given up to isomorphims, namely the algebra $\mathrm{Hom}(X,X)$
for whatever object $X$. We shall quantize $C(\hbar)$ into a stack
of linear categories, by first describing a prestack of algebroids
and then stackifying.

\subsection{The case of contractible foliations}

Let us suppose that the foliation $(N,F)$ is isomorphic to a fibration
of the type $M_{0}\times\mathbb{R}^{k}\rightarrow M_{0}$ for some $M_{0}$
and $k$. The point is that every foliation is locally of this form.
We shall construct an algebroid $A(N)$. The objects of this category
are the cross-sections%
\footnote{i.e.~submanifolds of $N$ intersecting each leaf transversally and
once (in our case they can be, of course, identified with maps $M_{0}\rightarrow\mathbb{R}^{k}$)%
} of the foliation. Morphisms are defined as follows (this construction,
together with the very definition of algebroids, is taken from \cite{K2}):
For any object $X$, the algebra $\mathrm{Hom}(X,X)$ is the quantized
algebra of functions on $X$. To define the space $\mathrm{Hom}(X,Y)$,
which should be a $\mathrm{Hom}(X,X)$-$\mathrm{Hom}(Y,Y)$ bimodule,
we first choose a smooth homotopy that moves $X$ to $Y$, while every
point stays in its leaf. This homotopy gives us an isomorphism between
the algebras $\mathrm{Hom}(X,X)$ and $\mathrm{Hom}(Y,Y)$, and we
define $\mathrm{Hom}(X,Y)$ to be the diagonal bimodule (i.e.~the
graph of the isomorphism). If we choose a different homotopy from
$X$ to $Y$, we have to identify the two definitions of $\mathrm{Hom}(X,Y)$.
To do it, we choose a homotopy between these homotopies, which gives
us an element $a\in\mathrm{Hom}(X,X)$ (see Section \ref{sec:Quantization-of-formal}),
and we identify the two $\mathrm{Hom}(X,Y)$'s by multiplication by
$a$. The element $a$ doesn't depend on the choice of homotopy of
homotopies, since any two such choices are homotopic.

\subsection{Construction of the stack}

Let us return to the case of a general $(N,F)$. For any open subset
$U\subset N$, on which the foliation satisfies the condition of the
previous subsection, we get the algebroid $A(U)$. Open sets of this
kind form a category $\mathit{Op_{res}}$: there is a (unique) morphism
$V\rightarrow U$ between two such sets iff $V\subset U$. We first define
a linear category $A$ fibred over $\mathit{Op_{res}}$, whose fibres
are $A(U)$'s: whenever $V\subset U$, for any object $X\in A(V)$
(i.e.~any cross-section of $V$) and $Y\in A(U)$, we define $\mathrm{Hom}(X,Y)$
by choosing a homotopy in $U$ that moves $X$ along the leaves to
a (uniquely defined) open subset $\bar Y$ of $Y$, and then continuing as in
the previous subsection, i.e.~setting $\mathrm{Hom}(X,Y)$ to be the graph of
the isomorphism between the quantized algebras of functions on 
$X$ and $\bar Y$.

The rest is purely formal. First we extend $A$ to a linear category
$A'$ fibred over $\mathit{Op}$ (the category of all open sets in
$N$). This is done by induction (right adjoint to the restriction).
Recall from \cite{gir} its possible construction: Any object $U'\in\mathit{Op}$
gives a full subcategory $\mathit{Op_{res}}(U')$ of $\mathit{Op_{res}}$,
of objects contained in $U'$. One then defines $A'(U')=\mathrm{Hom}_{/\mathit{Op_{res}}}(\mathit{Op_{res}}(U'),A)$.
We thus get a prestack over $N$ (a category fibred
over $\mathit{Op}$), which we finally stackify.

\end{document}